\documentclass[12pt, draft]{amsart}
\usepackage{amsfonts,amsmath,amsthm,latexsym,amscd,amsbsy,amssymb,euscript,fleqn,leqno,pb-diagram}%lamsarrow

\newcommand{\rdim}{rdim}

\newcommand{\diam}{diam}

\newcommand{\coker}{Coker}

\newtheorem{thm}{Theorem}[section]
\newtheorem{cor}[thm]{Corollary}
\newtheorem{lem}[thm]{Lemma}
\newtheorem{pro}[thm]{Proposition}

\theoremstyle{definition}

\theoremstyle{remark}

\newcommand{\ds}{\displaystyle}

\chardef\bslash=`\\ % p. 424, TeXbook

\makeatletter
\def\verbatim{\interlinepenalty\@M \@verbatim
  \leftskip\@totalleftmargin\advance\leftskip2pc
  \frenchspacing\@vobeyspaces \@xverbatim}
\makeatother
\numberwithin{equation}{section}

\begin{document}

%%%%%%%%%%% Begin Topmatter %%%%%%%%%%%%%%%%%

\title[Dimension-raising theorems for cohomological and extension dimensions]
{Dimension-raising theorems for cohomological and extension dimensions}

\author{Gencho Skordev}
\address{Center for Complex Systems and Visualization, Department of Computer Science and Mathematics,
University of Bremen, Universitatsalle 29, 28359 Bremen, Germany} \email{skordev@cevis.uni-bremen.de}

\author{Vesko Valov}
\address{Department of Computer Science and Mathematics, Nipissing University,
100 College Drive, P.O. Box 5002, North Bay, ON, P1B 8L7, Canada}
\email{veskov@nipissingu.ca}
\thanks{The second author was partially supported by NSERC Grant 261914-03.}



%%%%%%%%%% End topmatter %%%%%%%%%%%%%%%%%%%%%

\begin{abstract}
We establish cohomological and extension dimension versions of the Hurewicz dimension-raising theorem
\end{abstract}

\maketitle

\markboth{G.~Skordev and V.~Valov}{Dimension-raising theorems}

%%%%%%%%%%%%%%%%%%%%%%%%%%%%%%%%%%%%%%%%%%%%%%%%%%%%%%%%%%%%%%%%
%%%%%%%%%%%%%%%%%%%%%%%%%%%%%%%%%%%%%%%%%%%%%%%%%%%%%%%%%%%%%%%%

\section{Introduction}
Consider  a finite-to-one closed and surjective map $f\colon X\to Y$. The multiplicity function 
 $\mu\colon Y \to\mathbb N$ of the map $f$ is defined by $\mu (y)= card(f^{-1}(y))$, and let $\mu(f)= \{\mu (y): y \in Y\}$.
W. Hurewicz \cite{Hmu} proved   that $card(\mu(f)) \geq k+1$ provided 
$\dim Y \geq\dim X + k$ with $X$ and $Y$ being separable metric spaces. Moreover, if 
$ \mu (f) =\{ m(1)<m(2)<\ldots <m(\ell) \}$, where $\ell\geq k+1$, then 
$\displaystyle\dim Y_{m(s+1)}(f) \geq\dim Y -s$ for $ s=0, 1, \ldots , k$. Here
$Y_m(f) = \{y : \mu(y)=card (f^{-1}(y))\geq m\}$. The last assertion is a generalization of the following result  of H. Freudenthal \cite{F}: 
$\displaystyle\dim Y_{s+1}(f) \geq\dim Y -s$, $s=0,\ldots k $, if $\dim Y\geq\dim X + k$. Hence, 
$ Y_{k+1}(f)$ is not empty which was established first in \cite{H}.
The theorem of W. Hurewicz was further generalized for metric spaces in \cite{K, Sgod}, see also \cite{Filip, Szam, Ssb,Pears, E} for related results. 

The present paper deals with cohomological and extension dimension versions of the Hurewicz theorem.
Recall that the cohomological dimension $\dim_G$ of the paracompact space $X$ with coefficients in an Abelian group $G$ is defined as follows:
$\dim_G X\leq n$ if and only if the \v{C}ech cohomology $H^{n+1}(X, A, G) =0$ for every closed set $A$ in $X$. Equivalently, $\dim_G X \leq n$ if and only if the homomorphism $H^n(X, G)\to H^n(A, G)$ of the \v{C}ech cohomology groups induced by the inclusion $A \subset X$   is an epimorphism for every closed $A$ in $X$, see \cite{Kuzm, Dcohd, DArx}. 
If
$Z$ is a subset of a space $X$, we define $\rdim_{X,G}Z=\max\{\dim_GF: F\subset Z\hbox{}~\mbox{is closed in $X$}\}$.

There exists also another way to define $\dim_G$ for paracompact spaces. If $L$ is a $CW$-complex, we say that a space $X$ belongs to the class $\alpha(K)$ if every map $g\colon A\to K$, where $A\subset X$ is closed, can be extended to the whole space $X$ provided $g$ is extendable to some open neighborhood of $A$ in $X$. Shvedov \cite{sv} has shown that if $X$ is paracompact, then $\dim_GX\leq n$ if and only if $X\in\alpha(K(G,n))$, where $K(G,n)$ denotes the Eilenberg-MacLane complex of $G$ in dimension $n$ (for countable groups $G$ this was established by Huber \cite{hu}).

Concerning the definition of extension dimension for paracompact spaces, we adopt Dydak's approach \cite{dydak2}The : If $X$ is a paracompact space and $L$ a $CW$-complex, we write $e-\dim X\leq L$ if every map $f\colon A\to L$, where $A\subset X$ is closed, extends over $X$ up to homotopy. In such a case we say that $L$ is an absolute extensor of $X$ up to homotopy. Let us mention that, according to \cite[Proposition 3.5]{dydak2}, for every $CW$-complex $L$ there exists a complete metrizable simplicial complex $K(L)$ homotopy  equivalent to $L$ such that, for any paracompact space $X$ we have $K(L)$ is an absolute extensor of $X$ if and only if $L$ is an absolute extensor of $X$ up to homotopy.  
In case $X$ is compact or metrizable, the above definition of $e-\dim$ coincides with the original one introduced by Dranishnikov \cite{DExt},
i.e., every map $f\colon A\to L$ with $A\subset X$ being closed admits an extension over $X$. The notation $e-\dim X\leq e-\dim Y$ means that $e-\dim X\leq L$ for every $CW$-complex $L$ with $e-\dim Y\leq L$.

The paper is organized as follows. A cohomological version of the Hurewicz theorem is established in Section 2. Section 3 is devoted to a cohomological version of the mentioned above Freudenthal's theorem. In the final Section 4 we deal with dimension-raising theorems for extension dimension. The fist result in this section is an extension analogue of another dimension-raising theorem of Hurewicz stating that if $f\colon X\to Y$ is a closed map with $card(f^{-1}(y))\leq n+1$ for every $y\in Y$, then $\dim Y\leq\dim X+n$. This theorem, established first by Hurewicz \cite{H} for separable metric spaces,  has many generalizations. For the dimension $\dim$ and normal spaces it was proved by Zarelua \cite{ZSib}, Filippov \cite{Filip} and 
Pasynkov \cite{Pasynkov}.
Kuzminov \cite[Theorem 14, 1), p. 24]{Kuzm} was who first provided a cohomological version of this theorem for closed maps between finite dimensional metrizable compact spaces and cohomological dimension $dim_G$ with respect to an Abelian group $G$. Kuzminov's proof was based on test spaces for $dim_G$. With the same method he also established this theorem in the class of paracompact spaces with $G$ being either a periodic  group or a field \cite[p.39]{Kuzm}. Zarelua \cite{ZSib} introduced a new technique for investigating closed finite-to-one maps of paracompact spaces. As a byproduct, he obtained another generalization of the Hurewicz theorem for paracompact spaces and $dim_G$, where $G$ is a commutative ring with a unity.  Finally, for arbitrary Abelian groups $G$,  paracompact spaces, and more general maps, this theorem  was obtained by the first author \cite[Corollary 4.1]{Ssb}.

The introduction of extension dimension as a unification of the cohomological dimension $\dim_G$ and the ordinary dimension $\dim$ gave another way of possible generalizations of the Hurewicz results. Such a generalization was established by 
Dranishnikov and Uspenskij \cite[Theorem 1.6]{DU} for metrizable compacta. Our Proposition 4.1 and Corollary 4.2 extend the Dranishnikov-Uspenskij theorem for more general spaces. Let us mention that, in case of finite-dimensional compacta (resp., metrizable spaces) and simply connected $CW$-complexes, Proposition 4.1 and Corollary 4.2 follow from the first author's result mentioned above \cite[Corollary 4.1]{Ssb} combined with Dranishnikov's homological criterion for extensivity \cite[Theorem 9]{DExt} (resp., with Dydak's criterion \cite[Theorems H and G]{Dydak1}).  We also provide an extension version of the Hurewicz theorem
when $\mu(f)$ has finitely many values.

Everywhere below, by a group we mean an Abelian group. Recall that a map $f\colon X\to Y$ is perfect if $f$ is a closed  map having compact fibers $f^{-1}(y)$, $y\in Y$. It is well known that this is equivalent to $f$ being closed and $f^{-1}(K)$ compact for all compact sets $K\subset Y$. 

\section {Hurewicz theorem for the cohomological dimension}

Here, we shall prove a  theorem    for maps increasing the  cohomological dimensions $\dim _G $ with coefficients in a   group $G$.

\begin{pro}
If $f\colon X \to Y$ is a finite-to-one closed surjection between the metric space $X$ and $Y$ and  $\dim_G X +k \leq \dim_G Y<\infty$,  then $card(\mu (f))\geq k+1$. Moreover, if $\mu(f)=\{m(1) < m(2)< \ldots < m(\ell)\}$, then $\ell\geq k+1$, 
and $\displaystyle\dim_GY_{m(s+1)}(f) \geq\dim_G Y -s$ for $ s=0, 1, \ldots ,k$. 

If $\dim _G X < \dim_G Y = \infty$, then $card(\mu (f))\geq \infty$ and $\dim_G Y_m(f)=\infty$ for $ m=2,3, \ldots $.
\end{pro}

\begin{cor}
Let  $f\colon X \to Y$ be a finite-to-one closed  surjection between metrizable spaces and $card ( \mu (f))\leq k+1$. Then $\dim _G Y \leq\dim_G X + k $.
\end{cor}

\noindent
{\bf Remarks.} $1$. {\em Let  $f\colon X \to Y$ be a   closed  surjection between the metrizable spaces $X$ and $Y$. Then the
set $Y_m(f) = \{y : \mu(y)=card (f^{-1}(y))\geq m\}$ is a $F_ {\sigma}$ set in $Y$, see \cite[Lemma 4.3.5, p. 243]{E}. Therefore, $\dim_G Y_m(f) = \rdim_{Y,G}Y_m(f)$}.

\noindent
{\em $2$. Proposition $2.1$ and Corollary $2.2$ also hold in the more general situation when $X$ is paracompact and submetrizable  $($recall that a space is submetrizable if it admits a continuous metric, or equivalently, it admits a bijective continuous map onto a metrizable space$)$.  For the class of submetrizable spaces and the ordinary dimension $\dim$, Proposition $2.1$ and Corollary $2.2$ were established in \cite{Sgod} $($in this case the set $Y_m(f)$ is also
$F_\sigma$ in $Y$$)$. Bredon \cite[Theorem 8.15, p. 243]{B} also provided a version of Corollary $2.2$ for separable metric spaces $X, Y$ and a principal ideal domain $G$.}

The proof Proposition 2.1 is based on several lemmas. Before
starting the proof, let us provide some notions and results from cohomology and sheaf theory, see \cite{Godement, B, Sklyarenko}, and some constructions and theorems of A.V. Zarelua \cite{ZSib, ZTbil, ZKeld}. Recall that the cohomological dimension $\dim _{\mathcal L}Z$ of a space $Z$ with coefficients in a sheaf of Abelian groups ${\mathcal L}$ is 
defined by $\dim _{\mathcal L}Z = \min\{n : H^{n+1}(Z, {\mathcal L}_U)=0 \, \, \mbox{for every open set } \, \, U \subset Z \}$. Moreover,
$\dim _{\mathcal L}Z = \dim _G Z$ if ${\mathcal L}$ is the constant sheaf $Z\times G$.

First, we recall the notion of a local system of sheaves of groups on a given space $Y$, see \cite{ZSib, ZTbil, ZKeld}. Let $\Lambda$  be a partially ordered set. An open covering $\Omega = \{U_{\lambda}: \lambda \in \Lambda\}$  of $Y$ is said to be $Y$-directed \cite[Definition 7]{ZKeld} provided $U_{\mu} \subset U_{\lambda}$ for $ \lambda \leq \mu$ and  every index set $\Lambda _y = \{\lambda \in \Lambda : y \in U_{\lambda}\}$, $y\in Y$, is directed, i.e., if
$y\in U_{\lambda} \cap U_{\mu}$, then there exists $\nu \in \Lambda$ such that $\nu \geq \lambda , \, \nu \geq \mu$, and $y \in U_{\nu} \subset U_{\lambda} \cap U_{\mu}$. 

Let $\Omega = \{U_{\lambda}: \lambda \in \Lambda\}$ be an $Y$-directed open covering. A local system of sheaves 
$\Sigma = \{ \Omega , {\mathcal  L}_{\lambda}, \gamma_{\mu}^{\lambda}, \Lambda \}$ is a family of sheaves 
${\mathcal L}_{\lambda}$ on $U_{\lambda}$ and homomorphisms   
$\gamma_{\mu}^{\lambda}\colon {\mathcal L}_{\lambda}|U_{\mu} \to {\mathcal L}_{\mu}$ with $\mu \geq \lambda$   such that  $\nu  \geq \mu \geq   \lambda $ implies 
$\gamma ^{\lambda}_{\nu} = \gamma ^{\mu}_{\nu} \gamma ^{\lambda}_{\mu}$, see \cite[Definition 13]{ZKeld}.
A limit $limind \, \, \Sigma$ of the local system of sheaves $\Sigma$ is defined. This is a sheaf ${\mathcal L }$ on $Y$ with fibers 
${\mathcal L }_y = \lim ind \, \,  \{ ({\mathcal L_{\lambda}}) _y, \lambda \in \Lambda _y \}$ and the topology of ${\mathcal L }$ is induced by the topology of ${\mathcal L }_{\lambda}$, \cite[Definition 14]{ZKeld}.

We assume also that the homomorphisms $\gamma_{\mu}^{\lambda}$ are monomorphism, and that for every $U_{\lambda _{1}}, \, \, U_{\lambda _{2}}\in \Omega$ with $U_{\lambda _{1}} \cap U_{\lambda _{2}} \neq\emptyset $    there is $\mu \in \Lambda$ such that $ \mu \geq \lambda _1 , \lambda _2 $ and 
$U_{\mu}= U_{\lambda _{1}} \cap U_{\lambda _{2}}$. The local inductive systems of sheaves satisfying this condition are called regular, \cite[Definition 3.2]{ZTbil}.
 
 A local system of sheaves on $Y$, may be considered as a collection of inductive systems of sheaves parameterized by the points of the spaces $Y$. Then the limit of local system of sheaves may be interpreted as a collection of inductive limits of groups parameterized by the points of the space $Y$.

\begin{lem}\cite[Proposition 3.7]{ZSib}  Let $\Sigma = \{ \Omega , {\mathcal L }_{\lambda}, \gamma ^{\lambda}_{\mu}, \Lambda\}$ be a local   system of sheaves on the space $Y$. If $\rdim (Y, {\mathcal L }_{\lambda}) \leq n$, then $\dim _{{\mathcal L}}Y \leq n$, where ${\mathcal L} = lim ind \, \, \Sigma$.
\end{lem}
Here, $\rdim (Y, {\mathcal L }_{\lambda})=\max\{\dim_{\mathcal{L}_{\lambda}}F: F\subset U_\lambda \, \mbox{is a closed set in} \,  Y\}$. More generally,
for a fixed open subset $U$ of $Y$ and  a sheaf 
$\mathcal M$ on $U$, $\rdim (Y, {\mathcal M })$ denotes $\max \{\dim _{\mathcal M}  F : F \subset U , F \, \mbox{is  closed  in} \,  Y\}$. 

If $f\colon X\to Y$ is a closed surjective map and $G$ a group, then there is an exact sequence
$$0\to\mathcal G\to f_{\ast}f^{\ast}{\mathcal G} \to {\mathcal A} \to 0,\leqno{(1)}$$
where ${\mathcal G}= Y \times G$,    and  $f_{\ast}$, $f^{\ast}$ are the functors of direct and the inverse image of the sheaves with respect to the map $f$, see \cite{B}, Section I.3. Observe also that $f^{\ast}{\mathcal G} = X \times G$. 

Remind that a surjective map $f\colon X\to Y$ is called zero-dimensional if $\dim f^{-1}(y)=0$ for every $y \in Y$. This is equivalent to $\dim_G f^{-1}(y)=0$ for  every $y \in Y$ and some (arbitrary) group $G$. 

For closed zero-dimensional  maps
the sheaf ${\mathcal A}$ is concentrated on the set $Y_2(f)$, i.e.,  ${\mathcal A}_y=0$ for $y \notin Y_2(f)$.
This fact was used in \cite{Szam} and will be also exploited here. 

\begin{lem} Let $f\colon X \to Y$ be a perfect $0$-dimensional surjection on the paracompact space $X$. If $\dim_{G}X < \dim _G Y < \infty$, then
$\dim _{\mathcal A}Y = \dim _{G}Y -1$, where ${\mathcal A}$ is the sheaf from the exact sequence $(1)$. 
If $\dim_{G}X < \dim _G Y = \infty$, then $\dim _{\mathcal A}Y = \infty$.     
\end{lem}

\begin{proof}
Since the map $f$ is closed and $0$-dimensional, $f_{\ast}$ is an exact functor, see \cite[Proposition 1.4]{ZSib}.
Then, by \cite[Proposition 1.5]{ZSib}, 
$H ^i(X, {\mathcal G}_ {V})= H^i(Y, f_{\ast}f^{\ast}{\mathcal G}_U)$, where $V = f^{-1}(U)$,
${\mathcal G}_ {V} = (X \times G)_V=f^{\ast}{\mathcal G}_U$ and $U$ is an open set in $Y$. Therefore,
$\dim _{f_{\ast}f^{\ast}{\mathcal G}}Y =\dim_{G}X$. 

First, consider the case $\dim_GX<\dim_GY<\infty$.
The exact sequence (1) implies $\dim _{\mathcal A}Y \leq\dim_GY-1$.
If $m =\dim _GY$, we can find an open set $U$ in the space $Y$ such that  $H^m(Y, {\mathcal G}_U)\neq 0$. Then, 
according to (1), we have $H^{m-1}(Y, {\mathcal A}_U)\neq 0$ because $m>\dim_GY$. Hence,
$\dim _{\mathcal A} Y \geq m-1$.

Now  consider the case $\dim_G X < \dim_G Y = \infty$.
Here, $\dim_G Y = \infty$ means that there is a sequence of natural numbers $n_1 < n_2 < \ldots $ and open in $Y$ sets $U_k$ such that $H^{n_{k}}(Y, {\mathcal G}_{U_{k}}) \neq 0$ for $k=1, 2, \ldots$. Since 
$H^{n_{k}}(Y, f_{\ast}f^{\ast}{\mathcal G}_{U_{k}}) = H ^{n_{k}}(X, {\mathcal G}_ {V_{k}}) =0$ for $V_k = f^{-1}(U_k)$, 
the exact cohomology sequence induced by $(1)$ implies 
$$ 
H^{n_{k}-1}(Y, {\mathcal A}_{U_{k}}) \to  H^{n_{k}}(Y, {\mathcal G}_{U_{k}}) \to 0
$$ 
for every $n_k >\dim_G X$. Therefore, $H^{n_{k}-1}(Y, {\mathcal A}_{U_{k}}) \neq 0$ for all these $n_k$, which yields  $\dim _{\mathcal A}Y = \infty$. 
\end{proof}

Now, we describe a construction of A. V. Zarelua from his proof of \cite[Proposition 3.6]{ZTbil}. This construction is  given for closed finite-to-one and surjective maps $f\colon X\to Y$ of paracompact spaces $X$ and $Y$, but it also works for closed zero-dimensional surjective maps.

Local   systems 
$\Sigma _{\mathcal L},  \Sigma _{\mathcal M}, \Sigma _{ \tilde{\mathcal B}}$ of sheaves of groups on a space $Y$ are defined as follows. Let $\sigma = \{ U_1, \ldots U_{k(\sigma)} \}$ be a disjoint system of open sets in $X$ such that $f^{-1}(U) = U_1 \cup \ldots \cup U_{k(\sigma)} $ for some open set $U$ in $Y$. Consider the set $\Lambda = \{(U, \sigma) \}$ of all pairs $(U,\sigma)$ satisfying the above conditions and introduce a partial order on $\Lambda$ determined by the inclusion: $(U',\sigma')\geq (U,\sigma)$ if $U'\subset U$ and $\sigma'$ is a subdivision of $\sigma$.  Then  
$\Omega = \{U_{\lambda }= U:\lambda =(U, \sigma) \in \Lambda\}$ is an $Y$-directed open covering. 

Let $(U,\sigma)\in \Lambda , \, \sigma = \{U_1, \ldots U_{k(\sigma)}\}$ and 
$$
\Phi _j = \{ x \in U_j : f^{-1}f(x)\backslash U_j\neq\emptyset\},\, j=1,..,k(\sigma).
$$
Each set $\Phi _j$ is closed in $U_j$. Denote $F_j = f(\Phi_j)$ and let 
$F= F_1 \cup \ldots \cup F_{k(\sigma)}   $. Obviously, all sets $F$ and $F_j$, $j=1,\ldots,k(\sigma)$, are closed in $U$. Moreover, $F_j \subset Y_2(f), \, j= 1, \ldots , k(\sigma).$

Consider the sheaves  ${\mathcal L}(U, \sigma) = {\mathcal G}_F$, 
${\mathcal M}(U, \sigma) =  \oplus _1^{k(\sigma)}{\mathcal G}_{F_{j}}$ with ${\mathcal G} $ being the constant sheaf $U \times G$. The natural projections ${\mathcal G}_F \to {\mathcal G}_{F_{j}}$ induce a monomorphism 
$\tilde{\alpha}(U, \sigma): {\mathcal G}_F \to\oplus _1^{k(\sigma)}{\mathcal G}_{F_{j}}$.
The sheaf $\tilde{{\mathcal B}}(U, \sigma)$ is defined by the exact sequence
$$0 \to {\mathcal G}_F \to\oplus _{1}^{k(\sigma)}{\mathcal G}_{F_{j}} \to\tilde{\mathcal  B}(U, \sigma) \to 0,\leqno{(2)}$$ 
i.e., $\tilde{\mathcal  B}(U, \sigma) =\coker (\tilde{\alpha}(U, \sigma))$.

\smallskip
If $(U, \sigma), (U', \sigma ') \in \Lambda$ with $(U, \sigma) \leq  (U', \sigma ')$, there are natural homomorphisms 
$$
\gamma ^{(U, \sigma )}_{(U', \sigma '), {\mathcal L}}: {\mathcal L}(U, \sigma)|U'\to {\mathcal L}(U', \sigma ')  
$$
and
$$
\gamma ^{(U ,\sigma)}_{(U',\sigma'), {\mathcal M}}: {\mathcal M}(U, \sigma)|U' \to {\mathcal M}(U', \sigma ').
$$ 
The last two homomorphisms induce another one
$$
\gamma ^{(U,\sigma)}_{(U',\sigma'),\tilde{\mathcal B}}: \tilde{\mathcal B}(U,\sigma)|U' \to\tilde{\mathcal B}(U, \sigma).  
$$

Then 
$$
\Sigma _{\mathcal L}= \{ \Omega , {\mathcal L}(U, \sigma), \gamma ^{(U, \sigma )}_{(U', \sigma '),{\mathcal L}},\Lambda \},$$ 

$$\Sigma _{\mathcal M}= \{\Omega ,  {\mathcal M}(U, \sigma), \gamma ^{(U, \sigma )}_{(U', \sigma '), {\mathcal M}},\Lambda \}$$
and
$$\Sigma _{\tilde{\mathcal  B}}= \{ \tilde{\mathcal  B}(U, \sigma), \gamma ^{(U, \sigma )}_{(U', \sigma '), \tilde{\mathcal  B}},\Lambda \}
$$
are local inductive systems of sheaves on the space $Y$. 

Consider also the local inductive system 
$$\Sigma _{\tilde{\mathcal  D}}= \{ \tilde{\mathcal  D}(U, \sigma), \gamma ^{(U, \sigma )}_{(U', \sigma '), \tilde{\mathcal  D}}, \Lambda \},
$$
where $\tilde{\mathcal  D}(U, \sigma) = \oplus _1^{k(\sigma)}{\mathcal G}_{U{j}}$ and the homomorphisms 
$\gamma ^{(U, \sigma )}_{(U', \sigma '), \tilde{\mathcal  D}} $ are defined by natural projections.

Let ${\mathcal L}=limind \, \Sigma _{\mathcal L}$,  ${\mathcal M}=limind \, \Sigma _{\mathcal M}$ and 
$\tilde{\mathcal B}=limind \, \Sigma _{\tilde{\mathcal  B}}$.
 
\begin{lem}\cite[Proof of Proposition 3.6]{ZTbil}
Let $f: X \to Y$ be a perfect zero-dimensional surjection between paracompact spaces. Then the sheaves ${\mathcal A}$ and $\tilde{\mathcal  B}$ are isomorphic.
\end{lem}

\begin{proof}
In the case $f$ is finite-to-one, the lemma was established by A.V. Zarelua in \cite[Proof of Proposition 3.6]{ZTbil} after the construction of the local inductive systems $\Sigma _{\mathcal L}$, $\Sigma _{\mathcal M}$, $\Sigma _{\tilde{\mathcal B}}$ and their limits.  The crucial point in that proof, where the assumption that $f$ is finite-to-one is used,  is to show that the sheafs $ limind \, \Sigma _{\tilde{\mathcal  D}}$ and $f_{\ast}f^{\ast}{\mathcal G}$ are isomorphic. According to \cite[Proposition 1.4]{Ssb},
the last fact also holds when $f$ a zero-dimensional map with all fibers $f^{-1}(y)$, $y\in Y$, being  compact. The remaining part of the Zarelua arguments work in the present situation.
\end{proof} 

\begin{lem}
If $f: X \to Y$ be a perfect $0$-dimensional surjection and $X$ a paracompact space, then 
$\dim _{\mathcal A} Y \leq rdim_{Y, G} Y_2(f)$.
\end {lem}

\begin{proof}
The non-trivial case is $\rdim_{Y, G} Y_2(f) < \infty$.
According to Lemma 2.3, Lemma 2.5 and the exact sequence (2), it suffices to prove the inequalities 
$\rdim (Y, {\mathcal L}(U, \sigma)) \leq\rdim_{Y, G} Y_2(f)$ and 
$\rdim (Y, {\mathcal M}(U, \sigma)) \leq\rdim_{Y, G} Y_2(f)$ for 
any $(U,\sigma )\in \Lambda$. On the other hand, since $F= F_1 \cup \ldots \cup F_{k(\sigma)}$ and the sets $F$, \,    $F_j$, $j=1,\ldots,k(\sigma)$, are closed in $U$,
the second inequality implies the first one. So, we need to prove only the second inequality.

Let $C$ be a closed in $Y$ subset of $U$. Then 
$\dim _{{\mathcal G}_{F_{i}}}C = \dim _{{\mathcal G}_{F_i  \cap C}}C= \dim _G C \cap F_i$. Finally, since $C\cap F_i$ is a closed in $Y$ subset of $Y_2(f)$, we have
$\dim _G C \cap F_i\leq  rdim_{Y, G} Y_2(f)$. This completes the proof.
\end{proof}

\begin{cor}
Let $X$ be a paracompact space and $f: X \to Y$ a perfect $0$-dimensional surjection. If 
$\dim_G X<\dim _G Y<\infty$, then $\rdim _{Y, G} Y_2(f) \geq\dim _G Y -1.$ Moreover, if 
$\dim_G X<\dim _G Y=\infty$, then $\rdim _{Y, G} Y_2(f)=\infty$. 
\end{cor}
\begin{proof}
The proof follows from Lemma 2.4 and Lemma 2.6.
\end{proof}

\noindent
{\bf Remark.} {\em The second part of Corollary $2.7$ also holds when $f$ is a closed zero-dimensional surjection (not necessarily  perfect), but $X$ and $Y$ being metrizable.   
Indeed, then $\dim_{Y,G} Y_2(f)=\dim_G Y_2(f)$ and assume that $\dim_G Y_2(f) < \infty$. Let $Y^1 = Y \setminus Y_2(f)$ and $X^1 = f^{-1}Y^1$. The map $f_1 = f| X^1 : X^1 \to Y^1$ is a homeomorphism. Therefore, $\dim_G Y^1 = \dim_G X^1$. Moreover, $\dim_G X^1 \leq \dim_G X$, \cite[Theorem 16, 2]{Kuzm}. Hence, by  \cite[Theorem B]{dydak3},
$\dim_G Y \leq \dim_G Y^1 + \dim_G Y_2(f) + 2 < \infty$, which is a contradiction.}

\begin{lem}
Let $f: X \to Y$ be a closed $0$-dimensional surjection between metrizable spaces and $m\geq 2$. If 
$\dim_G X < \dim _G Y_m(f) < \infty$, then 
$\dim _{Y, G} Y_{m+1}(f) \geq\dim _{Y, G} Y_m(f) -1$. If
$\dim _{G}Y_m(f)=\infty$, then $\dim _{G}Y_{m+1}(f)=\infty$.
\end{lem}

\begin{proof}
The set $X_1 = \cup \{Bd(f^{-1}(y)) \, : \, y \in Y\}$, where $Bd(f^{-1}(y))$ is the boundary of $f^{-1}(y)$ in  $X$, is closed in $X$ because the map $f$ is closed, \cite[Lemma 23-2, p. 142]{N}. So, $f_1 = f| X_1 : X_1 \to Y$ is also closed,  surjective and zero-dimensional. By the Veinstein's lemma \cite[1.12.9, p. 111]{E}, all sets $f_1^{-1}(y)$,
$y \in Y$, are compact. Moreover,
$\dim_G X_1 \leq\dim_G X$. Therefore, without restriction of the generality, we can assume that the  original map $f$ 
is perfect.

First consider the case $\dim_G X < \dim_G Y_m(f) < \infty$.
Choose locally finite closed coverings $\omega^j = \{F_s^j\}_s$ of $X$ with 
$\diam F_s^j\leq 1/j$, $j= 1, 2,\ldots$. Since $Y_m(f)$ is an $F_{\sigma}$-set in $Y$, there is 
a closed set $C\subset Y$ with $C \subset Y_m(f)$ and $\dim _GC = \dim _G Y_m(f)$. Let $ \tilde{C}= f^{-1}(C)$
and $f_C = f|\tilde{C}: \tilde{C}\to C$. The map $f_C$ is perfect, zero-dimensional and surjective. Moreover,
$$Y_m(f_C)= \{y \in C: card(f_C^{-1}(y))\geq m\}= C.\leqno{(3)}$$
Let $y \in C$ and $\{x_1, \ldots , x_{m}\} \subset f^{-1}(y) = f^{-1}_C(y)$. There exists $j(y)$ such that
the covering $\omega ^{j(y)}$ separates the points $\{x_1, \ldots , x_{m}\}$, i.e., we can find
disjoint sets $F_{s_{1}}^{j(y)}, \ldots , F_{s_{m}}^{j(y)} \in \omega ^{j(y)}$ such that $x_i \in F_{s_{i}}^{j(y)}$, 
$i= 1, \ldots,m$. 

Consider the closed subsets $F_{j(y)}= \cap_{i=1}^{m} f_C(F_{s_{i}}^{j(y)} \cap \tilde{C})$ of $C$ and let
$C_j = \cup _{j(y)=j}F_{j(y)}$. Since $f$ is closed and the fibers $f^{-1}(y)$, $y\in Y$, are compact, all systems $\{ F_{j(y)}\}_{j(y)=j}$,
$y\in C$, are locally-finite and consist of closed sets. Therefore, each set $C_j$ is closed in $C$. 
Moreover, (3) implies $C = \cup _j C_j$. Then, there exist $j_0$ and $y_0 \in C$ with $j_0 = j(y_0)$ and
$\dim _{G} Y_m(f) = \dim _G C = \dim _G C_{j_{0}}= \dim _G F_{j(y_{0})}$. The sets
$F_{s_{1}}^{j_0}, \ldots , F_{s_{m}}^{j_0} \in \omega ^{j_0}$ separate the points $\{x_1, \ldots , x_{m}\}\subset f^{-1}(y_0)$. 
Let $\Phi = F_{s_{1}}^{j_0} \cap f^{-1}(F_{j(y_0)})$ and $f_{\Phi}= f| \Phi :\Phi \to F_{j(y_0)}$. The map $f_{\Phi}$ is perfect, surjective and $0$-dimensional and we have the following inclusion 
$$Y_2(f_\Phi) = \{ y \in F_{j(y_0)} : card(f_{\Phi}^{-1}(y))\geq 2 \} \subset Y_{m+1}(f).\leqno{(4)}$$
Moreover, 
$\dim _G \Phi\leq\dim _G X< \dim _{G} Y_m(f)=\dim _G F_{j(y_0)}<\infty$. Then, by Corollary 2.2, 
$\dim _{G} Y_2(f_\phi) \geq \dim_G F_{j(y_0)} -1=\dim _{G} Y_m(f) -1$. Finally, (4) completes the proof in the case $\dim_GY_m(f)<\infty$.

\smallskip
Now consider the case $dim_G Y_m(f)=\infty$.
Since $Y_m(f)$ is a $F_{\sigma}$ set in $Y$, we have two cases: (A)
There is a sequence of natural numbers $n_1 < n_2 < \ldots $ converging to $\infty$ and closed in $Y$ sets $C^k \subset Y_m(f)$ such that $\dim_G C^k = n_k$ for $k= 1, 2, \ldots$; (B)
There is a closed in $Y$ set $C \subset Y_m(f)$ such that $\dim_G C = \infty$.

{\em Case $($A$)$.}
Let $n_k > \dim_G X$. Applying the construction from the proof of the lemma in the case 
$\dim_G Y_(f) < \infty$ (with $C$ replaced by $C^k$), we obtain $\dim_G Y_{m+1}\geq n_k -1$. This implies $\dim_G Y_{m+1}= \infty$.

{\em Case $($B$)$.}
We apply again the construction from the case $\dim_G Y_m(f) < \infty $ for the set $C$ to obtain 
$C = \cup _j C_j$, where $C_j = \cup_{j(y)=j}F_{j(y)}$, $y\in C$. Since $\dim_G C = \infty$, we have again two cases:
$(B_1)$
There exists a sequence of natural numbers $n_1 < n_2 < \ldots $ converging to $\infty$ and closed sets $C_{m_k}\subset C$ such that $\dim_G C_{m_{k}} = n_k$; $(B_2)$
There exists $j_0$ such that $\dim_G C_{j_{0}}= \infty$ for some closed $C_{j_{0}}\subset C$.

In the   Case $(B_1)$  we argue as in the Case (A) and obtain $\dim_G Y_{m+1}= \infty$.

The Case $(B_2)$ splits again in two cases: $(B_{21})$
There is a sequence of natural numbers $n_1 < n_2 < \ldots$ such that $\lim n_k=\infty$ and $\dim_G F_{j(y_{k})}=n_k$ for $y_k \in C_{j_{0}}$ and $j(y_k) = j_0$; $(B_{22})$
There exists $y_0 \in C_{j_{0}}$ such that $\dim_G F_{j(y_0)}= \infty$.

The arguments from Case (A) applied under the hypotheses of  Case $(B_{21})$ imply $\dim_G Y_{m+1}= \infty$.

In Case $(B_{22})$ we apply the construction from the case $\dim_G Y_m(f) < \infty$ for the set $C$ and consider the map 
$f_{\Phi}: \Phi \to  F_{j(y_0)}$. This map is perfect, zero-dimensional and surjective. Moreover, $\dim_G \Phi \leq\dim_G X $ and $\dim_G F_{j(y_0)} = \infty$. Then, by  Corollary 2.7, $\dim_G Y_2(f_{\Phi}) = \infty$. Finally, inclusion $(4)$ yields $\dim_G Y_{m+1}= \infty$.
\end{proof}

\noindent
{\it Proof of Proposition $2.1$}.
Assume   $\dim _G X + k\leq\dim _G Y<\infty$, $k\geq 1$, and $\mu(f) = \{m_1 < \dots < m_{\ell}\}$. 
Then, by Corollary 2.7 and Lemma 2.8, 
$$\dim_{G} Y_{m_{k+1}}(f) \geq dim_G X.$$ 
So, $Y_{m_{k+1}}(f)$ is not empty, which yields $\ell\geq k+1$.

If $\dim_G X < \dim_G Y= \infty$, again by Corollary 2.7 and Lemma 2.8, we have $\dim_G Y_m(f) = \infty$ for $m \geq 2$. Therefore, $card (\mu(f))=\infty$.

\bigskip
Another result of Hurewicz \cite[Theorem II, p.74-76]{Hmu} is stating that if $f: X \to Y$ 
a closed surjective finite-to-one map between separable metrizable spaces, $X$ has the property 
$(\alpha)$ and $\dim Y - \dim X \geq n$,  then $card(\mu (f))\geq n+2$. 

Here, we say that a space $X$ has the property $(\alpha)$ if $\dim A <\dim X$ for every nowhere dense closed set $A$ in $X$. For example, every topological manifold $X$ satisfies the property $(\alpha)$.
In his proof Hurewicz used that the space $Y$ can be assumed to satisfy the following condition $(\beta)$: $\dim U = \dim Y$ for every open set $U$ in $Y$.

A version of this theorem also holds for the cohomological dimension. We say that $X \in (\alpha)_G$,
where $G$ is a group, if $\dim _{G}A <\dim_G X$ for every nowhere dense closed set $A\subset X$. Let us remind that every  locally compact n-cohomological manifold $M$ over the principal ideal domain $L$, notation $M$ is $n-cm_L$, possesses the property 
$(\alpha)_L$, see \cite[Proposition 4.9(a), p.14]{Borel}. 
We also consider the corresponding analogue of condition $(\beta)$: A space $Y$ satisfies condition $(\beta)_G$, notation $Y\in (\beta)_G$, if $\dim _G U =\dim _G Y$ for every open $U\subset Y$. 

Now, we can state a version of the Hurewicz result mentioned above. 

\begin{pro}
Let $f: X \to Y$ be a closed, surjective and finite-to-one map on the metrizable space $X$. If $X \in (\alpha)_G$, $\dim _G X + n \leq\dim _G Y < \infty$ and $Y \in (\beta )_G$, then $card(\mu (f))\geq n+2$.
\end{pro}

W. Hurewicz derived his theorem from a lemma whose cohomological version is Lemma 2.10 below. Observe that Proposition 2.9 follows from Lemma 2.10 and the inequality $(5)$.

\begin{lem}
Let $f: X \to Y$ be a closed finite-to-one map between the metrizable spaces $X$ and $Y$. Assume that $X \in (\alpha)_G$ and $Y\in(\beta)_G$. If $\dim_G X + k\leq\dim_G Y < \infty$ and  $m_0 = \max \{\mu (y) \, : y \in Y\}$, then $\dim_G Y_{m_{0}} < \dim _GX$.
\end{lem}

\begin{proof}
The proof of this lemma adapts  the arguments from the proof of the corresponding assertion for the dimension $\dim$, \cite[p.75-76]{Hmu}. Hurewicz worked with a countable basis of the space $X$. Instead of this, we  take the collection of elements of the closed coverings $\omega ^j$, see the proof of Lemma  2.8, Case (A). Then we can apply the constructions of W. Hurewicz with $\dim $ replaced by $\dim_G$.
\end{proof}

\section {Freudenthal's theorem for the cohomological dimension}
 In this section we shall prove the following
 
 \begin{pro}
 Let  $f: X \to Y$ be a perfect $0$-dimensional surjection between the paracompact spaces $X$ and $Y$. Suppose $\dim_GX +k\leq\dim_G Y < \infty $ for some natural number $k$ and a countable group $G$. Then 
 $\rdim_{Y, G} Y_{s+1}(f) \geq\dim_G Y - s$ for every $s=1,\ldots,k$.
 \end{pro}
 
 First, we need the following technical result. 
 
 \begin{lem}
 Let  $f: X\to Y$ be a $0$-dimensional closed and surjective map between the paracompact spaces $X$ and $Y$ and $G$ be a countable group. Suppose $F$ is a closed in $Y$ subset of $Y_2(f)$ and let $\Phi = f^{-1}(F).$ Then there exist families $\{\Phi_{\alpha, g}:\alpha\in\Lambda, g \in G\}$ and $\{F_{\alpha,g}:\alpha\in\Lambda, g\in G\}$ of closed sets in $X$ and $Y$, respectively, such that:
 \begin{enumerate}
 \item $\Phi = \cup _{\alpha , g}\Phi_{\alpha , g}$, $F = \cup _{\alpha , g} F_{\alpha , g}$ and 
 $f(\Phi _{\alpha , g})= F_{\alpha , g}$;

\smallskip
 \item $\{ F _{\alpha,g}:\alpha\in\Lambda, g \in G\}$ is locally countable in $F$ and, for every $g\in G$ the family $\{ F _{\alpha , g}:\alpha\in\Lambda\}$ is locally finite in $F$;

\smallskip
 \item $Y_{t}(\varphi _{\alpha , g}) = \{ y \in F _{\alpha , g} \, : \, card (\varphi _{\alpha , g}^{-1}(y)) \geq t \}\subset Y_{t+1}(f)$,
  where $\varphi _{\alpha , g}:\Phi _{\alpha , g}\to F_{\alpha , g}$ is the restriction map $f|\Phi_{\alpha,g}$.
 \end{enumerate}
 \end{lem}
 
 \begin{proof}
Recall that the support of the sheaf ${\mathcal A}$ is the set $Y_2(f)$. For every $y\in Y$ choose an open in $Y$ set $O_y$ and a section $s_y \in  \Gamma (\overline{O_y}, {\mathcal A})$ such that $s_y(z) \neq 0$ for all $z\in \overline{O_y} \cap F$. Here $\overline{O_y}$ is the closure of the set $O_y$ in $Y$. Assume that $O_y$ is so small that there exists a section $\tilde{s}_y \in \Gamma (\overline{O_y}, f_{\ast} G)$ such that $s_y$ corresponds to $\tilde{s}_y$ by the natural projection $f_{\ast} G \to {\mathcal A}$ from the exact sequence (1). Let $t_y \in \Gamma (f^{-1}O_y, G)$ be the section which corresponds to the section $\tilde{s}_y$ by the isomorphism of the groups of sections 
 $\Gamma (O_y , f_{\ast}G)$ and $\Gamma (f^{-1}O_y, G)$. Consider the covering $\omega = \{O_y \cap F\}_{y \in F}$ and 
 let $\tilde{\omega} = \{F_{\alpha}:\alpha\in\Lambda\}$ be a closed locally finite covering of $F$ refining $\omega$. Then, for every $\alpha$ there is $y(\alpha) \in F$ such that $F_ \alpha \subset \overline{O_{y(\alpha)}}$. Let $s_{\alpha} = s_{y(\alpha)}$, $\tilde{s}_{\alpha} = \tilde{s}_{y(\alpha)}$,  
 $t_{\alpha} = t_{y(\alpha)}$ and $\Phi _{\alpha}= f^{-1}(F_{\alpha})$. 

Consider the sets
 $\Phi _{\alpha , g}  = \{x \in \Phi _{\alpha} \, : \, t_{\alpha}(x)=g \}$, $\alpha\in\Lambda$, $g\in G$. Each  $\Phi _{\alpha , g}$ is closed in $X$ because the sheaf $X \times G$ is a Hausdorff space. Moreover, since $G$ is countable,  the family $\{ \Phi _{\alpha , g}:g\in G\}$ is countable for every $\alpha\in\Lambda$. Obviously, 
 $\Phi _{\alpha} = \cup _{g \in G} \Phi _{\alpha , g}$, and let $F_{\alpha , g} = f(\Phi _{\alpha , g})$. The family  $\{F_{\alpha , g}:\alpha\in\Lambda\}$ is locally finite in $F$ for every $g \in G$ 
because so is the family $\{F_{\alpha}:\alpha\in\Lambda\}$. This implies that $\{F_{\alpha , g}:\alpha\in\Lambda, g\in G\}$ is a locally countable family in $F$. Obviously, every element of the last family is a closed subset of $Y$. 
 
Denote by $\varphi _{\alpha , g}$ the restriction of $f$ on the set $\Phi _{\alpha , g}$. For every $y \in F_{\alpha , g}$ we have $s_{\alpha }(y) \neq 0$. Hence, there exist points $x_1, x_2 \in f^{-1}(y)$ such that 
 $t_{\alpha}(x_1) \neq t_{\alpha}(x_2)$. Indeed, otherwise $t_{\alpha}(x') = t_{\alpha}(x'')$ for all $x', x'' \in f^{-1}(y)$ would imply $s_{\alpha}(y)=0$, which is not possible. So, $x_i \in \Phi_{\alpha,g_i}$, where
 $g_i= t_{\alpha}(x_i)$, $i=1,2$. 
Therefore, 
 $$
 Y_{t}(\varphi _{\alpha , g}) = \{ y \in F _{\alpha , g} \, : \, card (\varphi _{\alpha , g}^{-1}(y)) \geq t \}\subset Y_{t+1}(f).   
 $$
This completes the proof.
\end{proof}

\noindent 
{\it Proof of Proposition $3.1$}

The proof goes by induction with respect to $k$. Corollary 2.7 provides the case $k=1$.
 Assume the assertion has been proved for all $k$ with $1 \leq k \leq\ell-1$ and let $k=\ell$.
 
 Choose a closed in $Y$ set $F$ which is contained in $Y_2(f)$ such that $\dim _G F =\rdim_{Y, G}Y_2(f)\geq\dim_GY -1 > dim_GX.$ Consider the families $\{\Phi_{\alpha, g}:\alpha\in\Lambda, g \in G\}$ and $\{F_{\alpha,g}:\alpha\in\Lambda, g\in G\}$ from Lemma 3.2. Since $F=\cup_{\alpha,g}F_{\alpha,g}$, $G$ is countable and $\{F_{\alpha,g}:\alpha\in\Lambda\}$ is locally finite in $F$ for every $g\in G$, $\dim _G F_{\alpha,g}=\dim _G F$ for some $\alpha\in\Lambda$ and $g\in G$.
Then both $F_{\alpha , g}$ and $\Phi _{\alpha , g}$ are paracompact spaces and the map 
 $f_{\alpha,g}: \Phi_{\alpha , g}\to F_{\alpha , g}$ is perfect, surjective and zero-dimensional. Moreover,
 $$
 \dim_G F_{\alpha , g}=\dim_G F \geq\dim_G Y -1 \geq\dim_GX + k -1\geq\dim _G \Phi _{\alpha , g} + k-1.
 $$
 It follows from the inductive assumption that
 $$
 \rdim_{Y, G}Y_s(f_{\alpha , g})=\rdim_{\Phi _{\alpha , g}, G} Y_{s}(f_{\alpha,g}) \geq\dim_G F_{\alpha , g} - s+1  
 $$
 for $s=1,\ldots, k-1$.
 According to Lemma 3.2(3), we have 
 $Y_s(f_{ \alpha , g}) \subset Y_{s+1}(f)$.
 Then
 $$
 \rdim_{Y, G} Y_{s+1}(f)\geq\rdim _{Y, G} Y_s(f_{\alpha , g})\geq\dim_G Y - s
 $$
 for $s=1,\ldots, k$,
 which is exactly what we need.

\section {Hurewicz theorem for extension dimension}
In this section we establish extension dimensional analogues of the Hurewicz theorems when either $card(f^{-1}(y))\leq n+1$ for every $y\in Y$ or $card(\mu(f))\leq n+1$. 
A map $f\colon X\to Y$ is said to have a metrizable kernel if there exist a metrizable space $M$ and a map $g\colon X\to M$ such that $g$ is injective on each fiber of $f$, i.e., $g|f^{-1}(y)$ is an one-to-one map for every $y\in Y$. Obviously, this is the case if $X$ is submetrizable. 

Proposition 4.1 below was established by Dranishnikov and Uspenskij \cite[Theorem 1.6]{DU} in the case both $X$ and $Y$ are metrizable compacta. Our proof of this proposition is a modification of the Dranishnikov and Uspenskij arguments.

\begin{pro}
Let $f\colon X\to Y$ be a closed surjective map admitting a metrizable kernel such that $Y$ is a paracompact $k$-space and $card(f^{-1}(y))\leq n+1$ for every $y\in Y$. Then $e-\dim Y\leq e-\dim (X\times\mathbb I^n)$.
\end{pro}

\begin{proof}
Let $m=n+1$ and $P_m(X)$ be the space of all probability measures on $X$ whose supports consist of at most $m$ points.
The map $f$ has an extension $P_m(f)\colon P_m(X)\to P_m(Y)$.  Actually, $P_m(f)$ is the restriction of the map
$P_m(\beta f)\colon P_m(\beta X)\to P_m(\beta Y)$ on $P_m(X)$, where $\beta f\colon\beta X\to \beta Y$ is the natural extension of $f$ between the
\v{C}ech-Stone compactifications of $X$ and $Y$. Since $f$ is a perfect map, so is $P_m(f)$. Hence,  
$Z=P_m(f)^{-1}(Y)$ is a paracompact $k$-space as a perfect preimage of $Y$. If $\triangle_n$ denotes the standard $n$-dimensional simplex in $\mathbb R^n$, there exists a map 
$p\colon X^m\times\triangle_n\to P_m(X)$ assigning to each point 
$(x_0,x_1,..,x_n,t_0,t_1,..,t_n)\in X^m\times\triangle_n$ the measure $\mu=\sum_{i=0}^{i=n}t_ix_i$. We consider  
another map $\pi\colon T\to X\times\triangle_n$, where $T=p^{-1}(Z)$, defined by $\pi(x_0,x_1,..,x_n,t_0,t_1,..,t_n)=(x_0,t_0,t_1,..,t_n)$. Obviously, $\pi$ is surjective.

\smallskip\noindent
{\em Claim. The maps $p|T$ and $\pi$ are perfect. Moreover $T$ is closed in $X^m\times\triangle_n$.}

Since $Z$ is a $k$-space, according to \cite[Theorem 3.7.18]{E1}, it suffices to show that $p^{-1}(K)$ is a compact subset of $T$ for every compact set $K\subset Z$. So, we fix such $K$ and consider the compact set  $F=f^{-1}(P_m(f)(K))\subset X$. Then $F^m\times\triangle_n$ is a compact subset of $X^m\times\triangle_n$ containing $p^{-1}(K)$. On the other hand
$p^{-1}(K)$ is closed in $X^m\times\triangle_n$ because $K$ is closed in $P_m(X)$. Therefore,  $p^{-1}(K)$ is compact.
Obviously, $X$ is $k$-space as a closed subset of $Z$, so is $X\times\triangle_n$. Consequently, to prove that $\pi$ is perfect,  again by \cite[Theorem 3.7.18]{E1}, it suffices to show that $\pi^{-1}(K\times\triangle_n)$ is compact for any compact $K\subset X$. Since 
$\pi^{-1}(K\times\triangle_n)$ is contained in $K^m\times\triangle_n$, the proof of the claim is reduced to establish that $T$ is closed in $X^m\times\triangle_n$. The last assertion follows from the fact that $Z$ is closed in $P_m(X)$ (because $Z=P_m(f)^{-1}(Y)$ and $Y$ is closed in $P_m(Y)$) and $T=p^{-1}(Z)$.

\smallskip
We are going to show that $e-\dim T\leq e-\dim(X\times\mathbb I^n)$. To this end, we need the following result of Dranishnikov and Uspenskij \cite[Theorem 1.2]{DU}: If $g\colon Z_1\to Z_2$ is a perfect $0$-dimensional surjection between paracompact spaces, then $e-\dim Z_1\leq e-\dim Z_2$ (this theorem was originally established for compact spaces, but the proof works for paracompact spaces as well). 
Since $\pi$ is perfect and, obviously, each fiber of $\pi$ is $0$-dimensional, the above theorem yields $e-\dim T\leq e-\dim(X\times\mathbb I^n)$. 

Next step is to show that $e-\dim Z\leq e-\dim(X\times\mathbb I^n)$.  We fix a metric space $(M,\rho)$ and a map $g\colon X\to M$ such that such $g$ is injective on each fiber $f^{-1}(y)$, $y\in Y$. Let $d$ be the continuous pseudometric on $X$ defined by $d(x',x'')=\rho(g(x'),g(x''))$ and consider the corresponding pseudometric $d_m$ on $X^m$, $d_m((x_i),(y_i))=\max_{0\leq i\leq n}d(x_i,y_i)$. For every equivalence relation $\Re$ on the set $\{0,1,..,n\}$ and every $\varepsilon>0$ let $E_{\Re,\varepsilon}\subset X^m$ be the set of all 
$m$-tuples $(x_0,..,x_n)$ satisfying the following conditions:

\begin{itemize}
\item if $i\Re j$, then $x_i=x_j$;
\item if $i$ is not $\Re$-equivalent to $j$, then $d(x_i,x_j)\geq\varepsilon$. 
\end{itemize}

\noindent
Since $d$ is a continuous pseudometric on $X$, $E_{\Re,\varepsilon}$ is a closed subset of $X^m$.
We choose one element from every equivalence class of $\Re$ and denote by $S_{\Re}$ the set of these $\Re$-representatives. Let also $$\triangle_{\Re}=\{t_0,..,t_n)\in\triangle_n: t_i>0\hbox{}~\mbox{if and only if}\hbox{}~i\in S_{\Re}\}$$  
As in the proof of Theorem 1.7 from \cite{DU}, one can show that the restriction of $p$ on the set $C\times\triangle_{\Re}$ is one-to-one
for every closed set $C\subset E_{\Re,\varepsilon}$ with a $d_m$-diameter $\leq\varepsilon/2$.  Since $\triangle_{\Re}$ is an $F_\sigma$-subset of $\triangle_n$, $(C\times\triangle_{\Re})\cap T$ is $F_\sigma$ in $T$. 

Now, for every $\varepsilon>0$ choose a locally finite closed covering $\omega_\varepsilon=\{H_{\alpha,\varepsilon}: \alpha\in\Lambda_\varepsilon\}$ of $X^m$ with each $H_{\alpha,\varepsilon}$ having a $d_m$-diameter $\leq\varepsilon/2$ 
(this can be done as follows: first, choose a locally finite closed covering $\gamma_{\varepsilon}$ of $M^m$ such that the $\rho_m$-diameter of each element of $\gamma_\varepsilon$ is $\leq\varepsilon/2$, where $\rho_m((a_i),(b_i))=\max_{0\leq i\leq n}\rho(a_i,b_i)$, and then let $\omega_\varepsilon$ to be $(g^m)^{-1}(\gamma_\varepsilon)$). Next, consider the locally finite in $T$ families $$\Theta_{\Re,k}=\{\big((H_{\alpha,1/k}\cap E_{\Re,1/k})\times\triangle_{\Re}\big)\cap T:\alpha\in\Lambda_{1/k}\}$$ with $k\in\mathbb N$ and $\Re$ being an equivalence relation on $\{0,1,..,n\}$. Each element of $\Theta_{\Re,k}$ is an $F_\sigma$-subset of $T$, so $\big((H_{\alpha,1/k}\cap E_{\Re,1/k})\times\triangle_{\Re}\big)\cap T=\bigcup_{j=1}^{\infty}H_{\Re,k,j}(\alpha)$ such that all $H_{\Re,k,j}(\alpha)$ are closed in $T$. Therefore, we obtain countably many families 
$$\Omega_{\Re,k,j}=\{p(H_{\Re,k,j}(\alpha)):\alpha\in\Lambda_{1/k}\}$$ 
of closed subsets of $Z$. Moreover, we already observed that $p$ restricted to $H_{\Re,k,j}(\alpha)$ is bijective. Because
$p|T$ is a closed map, this yields that $H_{\Re,k,j}(\alpha)$ and $p(H_{\Re,k,j}(\alpha))$ are homeomorphic. So, $e-\dim p(H_{\Re,k,j}(\alpha))\leq e-\dim(X\times\mathbb I^n)$. On the other hand, since all 
$\Theta_{\Re,k}$ are locally finite in $T$ and $p|T$ is perfect, the families $\Omega_{\Re,k,j}$ are locally finite in $Z$. Therefore, each set $Z_{\Re,k,j}=\cup\{p(H_{\Re,k,j}(\alpha)):\alpha\in\Lambda_{1/k}\}$ is closed in $Z$ and, according to the locally finite sum theorem for extension dimension \cite[Proposition 1.18]{DD}, $e-\dim Z_{\Re,k,j}\leq e-\dim(X\times\mathbb I^n)$ for any $\Re$ and $k,j\in\mathbb N$. Since $d$ is a metric on each fiber of $f$, it is easily seen that $Z$ is the union of the sets $Z_{\Re,k,j}$, $k,j\in\mathbb N$ and $\Re$ being an equivalence relation on $\{0,1,..,n\}$. Therefore, by the countable sum theorem for extension dimension, $e-\dim Z\leq e-\dim(X\times\mathbb I^n)$. 

The last step of our proof is to show that $e-\dim Y\leq e-\dim(X\times\mathbb I^n)$. We are going to use the following result of Dranishnikov and Uspenskij \cite[Proposition 2.3]{DU}:
Every surjective perfect map $g\colon Z_1\to Z_2$ with convex fibers between paracompact spaces has the following property: for any $CW$-complex $K$ and every closed subset $B\subset Z_2$ the restriction $g\colon g^{-1}(B)\to B$ induced a bijective map 
$g^*\colon [B,K]\to [g^{-1}(B),K]$ between the homotopy classes (such maps are called hereditary shape equivalences).
This result was established for compact spaces but its proof holds for paracompact spaces as well. Here, a map $g\colon Z_1\to Z_2$ has convex fibers if there exists a convex subset $E$ of a locally convex linear space and a closed embedding $j\colon Z_1\subset Z_2\times E$ such that the sets
$Z_1(y)=\{x\in E: (y,x)\in j(Z_1)\}$ are convex and compact for every $y\in Z_2$. The idea behind the proof of this result is the following simple fact: If $\phi\colon\bar{Y}\to\bar{Z}$ is an upper semicontinuous set-valued map with compact and convex values, where $\bar{Y}$ is paracompact and $\bar{Z}$ is a convex subset of a locally convex linear space, then for every family $\mathcal U$ of open in $\bar{Z}$ sets with each $\phi(y)$, $y\in\bar{Y}$, being a subset of some $U\in\mathcal U$, there exists a map $h\colon\bar{ Y}\to\bar{Z}$ such that  $h$ and $\phi$ are $\mathcal U$-close, i.e., for every $y\in\bar{Y}$ both $h(y)$ and $\phi(y)$ are contained in an element of $\mathcal U$. 

In our situation, we consider the perfect map $P_m(f)\colon Z\to Y$ which has convex fibers. Indeed, let $P(\beta X)$ be the space of all probability measures on $\beta X$. It is well know that $P(\beta X)$ is a compact and convex subset of $\mathbb R^{C(\beta X)}$, $C(\beta X)$ being the set of all continuous functions on $\beta X$. Moreover, there exists a natural map $P(\beta f)\colon P(\beta X)\to P(\beta Y)$ extending $P_m(f)$. It is easily seen that the map $j\colon Z\to Y\times P(\beta X)$, $j(\mu)=(P(\beta f)(\mu),\mu)$, is a closed embedding and all sets $P_m(f^{-1}(y))=\{\mu\in P(\beta X): (y,\mu)\in j(Z)\}$, $y\in Y$, are convex and compact. So, by the mentioned above result
\cite[Proposition 2.3]{DU}, $P_m(f)$ is a hereditary shape equivalence. Since $e-\dim Z\leq e-\dim(X\times\mathbb I^n)$,
this implies the required inequality $e-\dim Y\leq e-\dim(X\times\mathbb I^n)$.
\end{proof}

\begin{cor}
Let $f\colon X\to Y$ be a surjective map between the compact spaces $X$ and $Y$ with $card (f^{-1}(y) \leq n+1$ for every $y \in Y$. Then $e-\dim Y\leq e-\dim(X\times\mathbb I^n)$.  
\end{cor}

\begin{proof}
Suppose $L$ is a $CW$-complex such that $e-\dim(X\times\mathbb I^n)\leq L$ and $m=n+1$. Obviously, $X$ and $Y$ have the same topological weight. If, $w(X)=w(Y)$ is countable, the proof follows directly from Proposition 4.1. Otherwise, using the notations from the proof of Proposition 4.1, consider the spaces $Z=P_m(f)^{-1}(Y)$, $T=p^{-1}(Z)$ and the diagram 
$T\stackrel{p}{\longrightarrow}Z\stackrel{P_m(f)}{\longrightarrow}Y$. According to
the \v{S}\v{c}epin spectral theorem \cite{Sc}, we can find
continuous $\omega$-inverse systems $S_X=\{X_\alpha,\xi_{\alpha}^{\beta}: \alpha,\beta\in A\}$ and  
$S_Y=\{Y_\alpha,\varpi_{\alpha}^{\beta}: \alpha,\beta\in A\}$ consisting of metrizable compacta, and continuous maps
$f_\alpha\colon X_\alpha\to Y_\alpha$, $\alpha\in A$, such that $X=\underleftarrow{\lim}S_X$, $Y=\underleftarrow{\lim}S_Y$ and $\varpi_{\alpha}\circ f=
\xi_{\alpha}\circ f_\alpha$. Here, $\xi_{\alpha}\colon X\to X_\alpha$ and $\varpi_{\alpha}\colon Y\to Y_\alpha$ are the projections of the inverse systems $S_X$ and $S_Y$, respectively.
Then $P_m(X)=\underleftarrow{\lim}\{P_m(X_\alpha),P_m(\xi_{\alpha}^{\beta}): \alpha,\beta\in A\}$ and 
$X^m\times\triangle_n=\underleftarrow{\lim}\{{X_\alpha}^m\times\triangle_n, (\xi_{\alpha}^{\beta})^m\times id :\alpha,\beta\in A\}$. This implies 
$T=\underleftarrow{\lim}\{T_\alpha,(\xi_{\alpha}^{\beta})^m\times id: \alpha,\beta\in A\}$ and 
$Z=\underleftarrow{\lim}\{Z_\alpha,P_m(\xi_{\alpha}^{\beta}): \alpha,\beta\in A\}$, where
 $T_\alpha=(\xi_{\alpha}^m\times id)(T)$ and  $Z_\alpha=P_m(\xi_{\alpha})(Z)$. We also consider the maps $p_\alpha\colon T_\alpha\to Z_\alpha$ assigning to each point
$\big(\xi_{\alpha}(x_0),\xi_{\alpha}(x_1),..,\xi_{\alpha}(x_n),t_0,t_1,..,t_n\big)\in X_\alpha^m\times\triangle_n$ the measure 
$\sum_{i=0}^{i=n}t_i\xi_{\alpha}(x_i)\in Z_\alpha$.
Therefore,
the following diagrams are commutative for all
$\alpha$:

$$
\begin{CD}
T@>{p}>>Z@>{P_m(f)}>>Y@<{f}<<X@>{id}>> P_m(X)\\ @ VV{\xi_\alpha^m\times id}V
@VV{P_m(\xi_\alpha)}V @VV{\varpi_\alpha}V @VV{\xi_\alpha}V @VV{P_m(\xi_\alpha)}V\\
T_\alpha@>{p_\alpha}>>Z_\alpha@>{P_m(f_\alpha)}>>Y_\alpha@<{f_\alpha}<<X_\alpha @>{id}>>P_m(X_\alpha)
\end{CD}
$$

Let us observe that the requirement for $f$ in Proposition 4.1 to have a metrizable kernel was necessary only to prove that $e-\dim Z\leq e-\dim T$, the other two facts we established there remain valid in the present situation. So, 
$e-\dim T\leq e-\dim(X\times\mathbb I^n)$ and $e-\dim Y\leq e-\dim Z$. Thus, we need to show that 
$e-\dim Z\leq e-\dim(X\times\mathbb I^n)$.

Since $e-\dim T\leq e-\dim(X\times\mathbb I^n)\leq L$, we can use the factorization theorem for extension dimension (see \cite[Theorem 2]{pasynkov2}, \cite{lrs}) to construct the inverse system $S_T=\{Z_\alpha,P_m(\xi_{\alpha}^{\beta}): \alpha,\beta\in A\}$ in such a way that $e-\dim T_\alpha\leq L$, $\alpha\in A$. The equalities
$Z_\alpha= P_m(f_\alpha)^{-1}(Y_\alpha)$ and $T_\alpha = p_\alpha^{-1}(Z_\alpha)$ may not be true
 for some $\alpha$, but we always have the inclusions
$Z_\alpha\subset P_m(f_\alpha)^{-1}(Y_\alpha)$ and $T_\alpha\subset p_\alpha^{-1}(Z_\alpha)$. Since each $f_\alpha$ has a metrizable kernel, according to the arguments from the proof of Proposition 4.1, $e-\dim Z_\alpha\leq e-\dim T_\alpha$.  Consequently, $e-\dim Z_\alpha\leq L$ for every $\alpha$. Finally, $Z$ being the limit of the inverse system $S_T$ yields $e-\dim Z\leq L$.                                 
\end{proof}

Next proposition provides the relation between $e-\dim Y$ and $e-\dim X$ when $f\colon X\to Y$ is a finite-to-one closed map with $card(\mu(f))\leq n+1$. For any $CW$-complexes $K$ and $L$ their join is denoted by $K\star L$.

\begin{pro} Let $f\colon X\to Y$ be a closed surjection between metrizable spaces and $card(\mu(f))\leq n+1$. 
If $e-\dim X\leq L$, then $e-\dim Y\leq\underbrace{ L\star L\star...\star L}_{n+1}$.
\end{pro}

\begin{proof}
Let $\ds\mu(f)=\{k_1,k_2,..,k_{\ell}\}$ and $Y^i=\{y\in Y:card(f^{-1}(y)=k_i\}$ for all $i=1,..,\ell$. Obviously, 
$Y=\bigcup_{i=1}^{i=\ell}Y^{i}$.
Since the extension dimension of a union of two subsets of any metrizable space is $\leq$ the join of their extension dimensions \cite[Theorem A]{dydak3}, it suffices to prove that $e-\dim Y^{i}\leq L$ 
for any $i$.

Passing to the subsets $Y^i\subset Y$ and $X^i=f^{-1}(Y^i)\subset X$, we can assume that all fibers $f^{-1}(y)$, $y\in Y$, have cardinality exactly $m$ for some $m\geq 1$. Now, fix a metric $d$ on $X$ and consider the space $exp_m(X)=\{F\subset X: card(F)\leq m\}$ equipped with the Vietoris topology, or equivalently,
with the Hausdorff metric generated by $d$. There exists a natural continuous map $exp_m(f)\colon exp_m(X)\to exp_m(Y)$ defined by $exp_m(f)(F)=f(F)$. Obviously, $exp_m(f)$ is the restriction of the map 
$exp_m(\beta f)\colon exp_m(\beta X)\to exp_m(\beta Y)$ on $exp_m(X)$. Since $f$ is perfect, the last observation implies that $exp_m(f)$ is also perfect. So, $Z=exp_m(f)^{-1}(Y)$ is a closed subset of $exp_m(X)$ because $Y$ is closed in $exp_m(Y)$. 

There is also a map $p\colon X^m\to exp_m(X)$, $p((x_1,..,x_{m}))=\{x_1,..,x_{m}\}$. It is easily seen that $p$ is perfect, so is $p|T$, where  $T=p^{-1}(Z)$. Moreover, $Z$ being closed in $exp_m(X)$ yields that $T$ is closed in $X^m$. Let us
consider the projection $\pi\colon T\to X$, $\pi((x_1,..,x_m))=x_1$. Since $(x_1,..,x_m)\in T$
if and only if $f(x_1)=...=f(x_m)$, we have $\pi^{-1}(K)\subset f^{-1}(f(K))^m\cap T$ for every $K\subset X$. This implies that $\pi$ is finite-to-one and $\pi^{-1}(K)$ is compact provided $K\subset X$ is compact. Hence, $\pi$ is a perfect
light map, and by \cite[Theorem 1.2]{DU}, $e-\dim T\leq e-\dim X\leq L$. 

Next step is to show that $e-\dim U\leq L$, where $U$ is the subset of $Z$ such that every $F\in U$ consists of $m$ points. Since $U$ is open in $Z$, so is $T_U=p^{-1}(U)$ in $T$. Consequently, by the countable sum theorem, 
$e-\dim T_U\leq e-\dim T\leq L$.
Now, for every $\varepsilon>0$  
let $$E_\varepsilon=\{(x_1,..,x_m)\in T_U: d(x_i,x_j)\geq\varepsilon\hbox{}~\text{for}\hbox{}~i\neq j\}.$$
If $C$ is a closed subset of $E_\varepsilon$ of $d_m$-diameter $\leq\varepsilon/2$, where $d_m$ is the metric on $X^m$ defined by the same equality as in Proposition 4.1,
then $p|C$ is a homeomorphism. Indeed, suppose $p((x_1,..,x_m))=p((y_1,..,y_m))$ for some $(x_1,..,x_m), (y_1,..,y_m)\in C$.
Then $\{x_1,..,x_m\}=\{y_1,..,y_m\}$ and, since $d(x_i,y_i)\leq\varepsilon/2$ for all $i$, it follows that
$x_i=y_i$, $i\in\{1,..,m\}$. So, $p|C$ is injective. On the other hand $p\colon T_U\to U$ is a perfect map and $C$ is closed in $T_U$ because $E_\varepsilon\subset T_U$ is closed. Hence, $p|C$ is a homeomorphism with $p(C)$ being a closed
subset of $U$. Obviously, the sets $E_{1/k}$, $k\geq 1$, form a closed covering of $T_U$. For each $k$ consider
a locally finite in $T_U$ covering $\Theta_k$ of $E_{1/k}$ with each $H\in\Theta_k$ being a closed set in $T_U$ of $d_m$-diameter $\leq 1/2k$. Therefore, $\{p(H):H\in\Theta_k\}$ is a locally finite family in $U$ covering $p(E_{1/k})$ and consisting of closed sets in $U$ with $e-\dim p(H)\leq L$. Consequently, $D_k=\bigcup\{p(H):H\in\Theta_k\}$ is closed in $U$ and
$e-\dim D_k\leq L$. Finally, since $U=\cup_{k=1}^{\infty}D_k$, we obtain $e-\dim U\leq L$. 
  
The last step is to show that $e-\dim Y\leq L$. We define the sets 
$$Y_\varepsilon=\{y\in Y:d(x',x'')>\varepsilon\hbox{}~\text{for every}\hbox{}~x'\neq x''\in f^{-1}(y)\},$$ $\varepsilon>0$. Using that $f$ is a perfect map and any fibers of $f$ contains exactly $m$ points, one can show that 
each $Y_\varepsilon$ is open in $Y$ and all maps $g_\varepsilon\colon Y_\varepsilon\to U$, $g_\varepsilon(y)=f^{-1}(y)$, 
are embeddings. Hence, $e-\dim Y_\varepsilon\leq L$ because $e-\dim U\leq L$ yields $e-\dim A\leq L$ for every subset $A\subset U$. To complete the proof, observe that $Y=\bigcup_{k=1}^{\infty}Y_{1/k}$ with each $Y_{1/k}$ being an $F_\sigma$ subset of $Y$ implies $e-\dim Y\leq L$.
\end{proof}

{\bf Acknowledgments.} The authors would like to express their sincere thanks to E.G.Sklyarenko and A.V.Zarelua for fruitful and stimulating discussions.

\end{document}